%% file: nothomotohyp.tex
\documentclass[11pt,notitlepage, leqno]{amsart}
  
 \include{packages}

 \renewcommand{\tilde}{\widetilde}

\usepackage{subfig}
\usepackage{floatrow}

\makeatletter 
\def\paragraph{\@startsection{paragraph}{4}%
\z@\z@{-\fontdimen2\font}%
{\normalfont\bfseries}}
\makeatother

\begin{document} 
	  	\title[]{A locally hyperbolic 3-manifold that is not homotopy equivalent to any hyperbolic 3-manifold}
	 	\author{Tommaso Cremaschi}
	 	\date{\today}
 		\thanks{The author gratefully acknowledges support from the U.S. National Science Foundation grant DMS-1564410: Geometric Structures on Higher Teichm\"uller Spaces.}

	 	\maketitle
	 	
		\small
	 	
		\paragraph{Abstract:} We construct a locally hyperbolic 3-manifold $M$ such that $\pi_ 1(M)$ has no divisible subgroups. We then show that $M$ is not homotopy equivalent to any complete hyperbolic manifold.
		\normalsize
	 \section*{Introduction}
	 
	 In \cite{C2017} we constructed an example of a locally hyperbolic 3-manifold without divisible elements that was not homeomorphic to any complete hyperbolic 3-manifold. This answered a question of Agol \cite{DHM,Ma2007}. In the present work we show how, using similar techniques to \cite{C2017}, one can construct a locally hyperbolic 3-manifold $M$ without divisible elements such that $M$ is not homotopy equivalent to any complete hyperbolic 3-manifold:
	 \begin{Theorem*}\label{maintheorem} There exists a 3-manifold $M$ that is locally hyperbolic without divisible elements in $\pi_1(M)$ that is not homotopy equivalent to any hyperbolic 3-manifold.
\end{Theorem*}

	\paragraph*{Acknowledgements:}  I would like to thank Ian Biringer and Martin Bridgeman for many helpful discussions.

	 	\vspace{0.5cm}

		\paragraph*{Notation:} We use $\simeq$ for homotopic and by $\pi_0(X)$ we intend the connected components of $X$. With $\Sigma_{g,k}$ we denote the genus $g$ orientable surface with $k$ boundary components. By $N\hookrightarrow M$ we denote embeddings while $S\imm M$ denotes immersions.
		
		\vspace{0.5cm}

		\section{Background} We now recall some facts and definitions about the topology of 3-manifolds, more details can be found in \cite{He1976,Ha2007,Ja1980}.

An orientable 3-manifold $M$ is said to be \emph{irreducible} if every embedded sphere $\mathbb S^2$ bounds a 3-ball. A map between manifolds is said to be \emph{proper} if it sends boundaries to boundaries and pre-images of compact sets are compact. We say that a connected properly immersed surface $S\imm M$ is \emph{$\pi_1$-injective} if the induced map on the fundamental groups is injective. Furthermore, if $S\hookrightarrow M$ is embedded and $\pi_1$-injective we say that it is \emph{incompressible}. If $S\hookrightarrow M$ is a non $\pi_1$-injective two-sided\footnote{The normal bundle on $S$ in $M$ is trivial.} surface by the Loop Theorem we have that there is a compressing disk $D\hookrightarrow M$ such that $\partial D=D\cap S$ and $\partial D$ is non-trivial in $\pi_1(S)$. 

An irreducible 3-manifold $(M,\partial M)$ is said to have \emph{incompressible} \emph{boundary} if every map: $(D^2,\partial D^2)\hookrightarrow (M,\partial M)$ is homotopic via a map of pairs into $\partial M$. Therefore, $(M,\partial M)$ has incompressible boundary if and only if each component $S\in\pi_0(S)$ is incompressible, that is $\pi_1$-injective.  An orientable, irreducible and compact $3$-manifold is called \textit{Haken} if it contains a two-sided $\pi_1$-injective surface. A 3-manifold is said to be \emph{acylindrical} if every map $(S^1\times I,\partial (S^1\times I))\rar (M,\partial M)$ can be homotoped into the boundary via maps of pairs.

\bdefi A 3-manifold $M$ is said to be \emph{tame} if it is homeomorphic to the interior of a compact 3-manifold $\overline M$. \edefi


\bdefi We say that a codimension zero submanifold $N\overset{\iota}{\hookrightarrow} M$ forms a \emph{Scott core} if the inclusion map $\iota$ is a homotopy equivalence. \edefi

		By \cite{Sc1973,HS1996,RS1990} given an orientable irreducible 3-manifold $M$ with finitely generated fundamental group a Scott core exists and is unique up to homeomorphism.		
		
Let $M$ be a tame 3-manifold, then given a Scott core $C\hookrightarrow M\subset \overline M$ with incompressible boundary we have that, by Waldhausen's cobordism Theorem \cite{Wa1968}, every component of $\overline{\overline M\setminus C}$ is a product submanifold homeomorphic to $S\times I$ for $S\in\pi_0(\partial C)$.

		\bdefi Given a core $C\hookrightarrow M$ we say that an end $E\subset\overline{ M\setminus C}$ is \textit{tame} if it is homeomorphic to $S\times [0,\infty)$ for $S=\partial E$.\edefi
		
		 A core $C\subset M$ gives us a bijective correspondence between the ends of $M$ and the components of $\partial C$. We say that a surface $S\in\pi_0(\partial C)$ \textit{faces the end $E$} if $E$ is the component of $\overline{M\setminus C}$ with boundary $S$. It is a simple observation that if an end $E$ facing $S$ is exhausted by submanifolds homeomorphic to $S\times I$ then it is a tame end.

\bdefi
A 3-manifold $M$ is said to be \emph{locally hyperbolic} if every cover $N\twoheadrightarrow M$ with $\pi_1(N)$ finitely generated is hyperbolizable.
\edefi
\brem
By the Tameness and Geometrization Theorem \cite{Kap2001,Per2003.1,Per2003.2,Per2003.3} this is equivalent to saying that every cover $N\twoheadrightarrow M$ with $\pi_1(N)$ finitely generated is the interior of a compact, irreducible 3-manifold that is atoroidal and with infinite $\pi_1$.
\erem
\bdefi An element $g$ in a group $G$ is said to be \emph{divisible} if for all $n\in\N$ there exists $\alpha\in G\setminus \set{e} $ such that $g=\alpha^n$.
\edefi
\subsection{Homotopy equivalences}
We now recall some facts about homotopy equivalences of irreducible 3-manifolds. For details see \cite{Jo1979,JS1978,CM2006}.

\bdefi
 A \emph{Seifert fibered 3-manifold $M$} is a compact, orientable, irreducible 3-manifold that has a fibration by circles.
\edefi

\bdefi
Given a 3-manifold $M$ and a mnaifold $N$, $\dim (N)\leq 3$, a continuous map $f:(N,\partial N)\rar (M,\partial M)$ is \emph{essential} if $f_*$ is $\pi_1$-injective and $f$ is not homotopic via map of pairs to a map $g$ such that $g(T)\subset\partial M$.
\edefi

\bdefi
Given an irreducible, compact 3-manifold $(M,\partial M)$ with incompressible boundary a \emph{characteristic submanifold} for $M$ is a codimension zero submanifold $(N,R)\hookrightarrow ( M,\partial  M)$ satisfying the following properties:
\begin{enumerate}
\item[(i)] every component $(\Sigma,\partial \Sigma)\in\pi_0(N)$ is an $I$-bundle or a Seifert fibered manifold;
\item[(ii)] $\partial N\cap \partial M=R$;
\item[(iii)] all essential maps of a Seifert fibered manifold $S$ into $ ( M,\partial  M)$ are homotopic as maps of pairs into $(N,R)$.

\end{enumerate}

\edefi
By work of Johannson and Jaco-Shalen  we have such a submanifold for compact, irreducible 3-manifolds with incompressible boundary \cite[2.9.1]{CM2006}:

\begin{Theorem*}[Existence and Uniqueness] Let $(M,\partial M)$ be a compact, irreducible 3-manifold with incompressible boundary. Then there exists a characteristic submanifold $(N,R)\hookrightarrow (M,\partial M)$ and any two characteristic submanifolds are isotopic. 
\end{Theorem*}

This is also called the JSJ or annulus-torus decomposition. One of the main application of the JSJ decomposition is the following Theorem \cite[24.2]{Jo1979}:

\bthm\label{homequivclass}
Let $(M,\partial M)$ and $(M',\partial M')$ be compact irreducible 3-manifolds with incompressible boundary and denote by $(N,R)$, $(N',R')$ respectively their characteristic submanifolds. Given a homotopy equivalence $f:M\rar M'$ then we have a map $\phi$ homotopic to $f$ such that:
\begin{enumerate}
\item[(i)] $\phi:\overline{M\setminus N}\diffeo \overline{M'\setminus N'}$ is a homeomorphism;
\item[(ii)] $\phi:N\rar N'$ is a homotopy equivalence.
\end{enumerate}
\ethm

In particular if $M$ is acylindrical we have that $M$ has no characteristic submanifold therefore any homotopy equivalence $f:M\rar N$ is homotopic to a homeomorphism.

\bdefi
Given an essential properly embedded annulus $(A,\partial A)$ in $(M,\partial M)$ a \emph{Dehn flip} of $M$ along $A$ is the 3-manifold $N$ obtained by cutting $M$ along $A$ picking a homeomorphism $f:M\vert A\diffeo M\vert A$ that is the identity on $A$ and re-gluing $f(M\vert A)$ along $f(A)$ either via the identity or via the map $\phi(x,t)=(x, 1-t)$ where we parametrised $A$ by $S^1\times [0,1]\cong A$.

A Dehn flip of $M$ along $A$ naturally gives a homotopy equivalence $h:M\rar N$ which we will also denote by a \emph{Dehn flip}.

\edefi

Moreover, as a consequence of Theorem \ref{homequivclass}, see\cite{Jo1979,JS1978,CM2006}, we have that homotopy equivalences of $M$ are generated by \emph{Dehn flips} along annuli contained in the boundary of thecharacteristic submanifold of $M$, see Theorem \cite[29.1]{Jo1979}. As a consequence of Johansson homotopy equivalence theory for Haken 3-manifold we get:

\blem\label{dehnflip} If the characteristic submanifold of a Haken 3-manifold $M$ is given by one embedded separating cylinder $C$ then any 3-manifold $N$ homotopy equivalent to $M$ is either homeomorphic to $M$ or to a Dehn flip along $C$.
\elem

	  \section{Construction of the Example}

	 Consider the 3-manifold $A$ obtained as a thickening of the 2-complex given by gluing a genus two surface $S$ and a torus $T$ so that a meridian of $T$ is identified with a separating simple closed curve $\gamma$ of $S$. Note that $\partial A$ is formed by two genus two surfaces both of which are incompressible in $A$. Let $B,C$ be two copies of a hyperbolizable, acylindrical 3-manifold with incompressible genus two boundary (for example see \cite[3.3.12]{Th1997}) and glue $B,C$ to the 3-manifold $A$ one to each boundary component. Then we obtain a closed 3-manifold $X$:

\begin{center}\begin{figure}[h!]
						\def\svgwidth{300pt}
						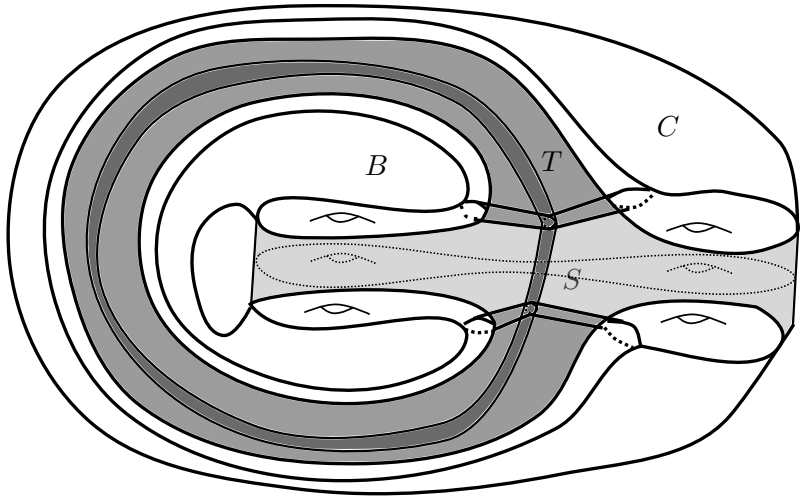
						\caption{Schematic of the manifold $X$ where the starting manifold $A$ is shaded and the essential torus $T$ is dark grey.}	
						\end{figure}
\end{center}
Note that the manifold $X$ is not hyperbolizable since it contains the essential torus $T$ and that the surface $S$ is incompressible and separating in $X$. 
\brem\label{fact0}
The 3-manifold $Y\eqdef X\vert S= X\setminus N_\epsilon(S)$ is hyperbolizable, with incompressible boundary and its characteristic submanifold is given by an annulus connecting the two distinct boundaries. Thus, any annulus with both boundary components on the same surface is boundary parallel.
\erem
The infinite cyclic cover $M$ of $X$ is obtained by gluing infinitely many copies $\set{Y_i}_{i\in\N}$ of $Y$ along their boundaries:

\begin{center}\begin{figure}[h!]
						\def\svgwidth{300pt}
						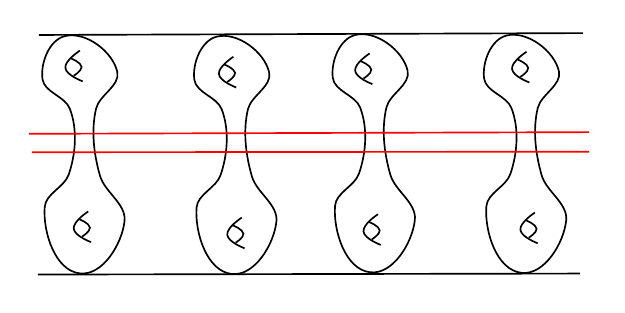
						\caption{}\label{fig1}
						\end{figure}
\end{center}

Denote by $\set{S_i}_{i\in\Z}\subset M$ the lifts of $S$. The surfaces $\set{S_i}_{i\in\Z}$ are all genus two and incompressible in $M$. Moreover, we denote by $M_{i,j}\eqdef \cup_{i\leq k\leq j} Y_k$ the compact submanifolds of $M$ co-bounded by $S_i,S_j$ for $i<j$ and by $A$ the properly embedded annulus in $M$ that is the lift of the essential torus $T$ and we let $A_{i,j}\eqdef M_{i,j}\cap A$. With an abuse of notation we denote by $\gamma\in\pi_1(M_{i,j})$ the elements corresponding to $\pi_1(A)$.

In the remainder of this work we will show that the manifold $M$ satisfies the following three properties:
\begin{enumerate}
\item $\pi_1(M)$ has no divisible elements;
\item $M$ is locally hyperbolic;
\item $M$ is not homotopy equivalent to any hyperbolic 3-manifold.
\end{enumerate}

Which will give us the theorem \ref{maintheorem}. We will now show that $M$ is locally hyperbolic and that $\pi_1(M)$ has no divisible elements.

\blem\label{nodivis} The manifold $M$ has no divisible elements in $\pi_1(M)$.
\elem

\bpf The manifold $M $ is the cover of a compact 3-manifold $X$ thus we have that $\pi_1(M)\subset \pi_1(X)$. Since $X$ is irreducible, compact and with infinite $\pi_1$ by \cite{Sh1975} we have that $\pi_1(X)$ it has no divisible elements in $\pi_1$. \epf

\blem\label{lochyp} The manifold $M$ is locally hyperbolic and all covers corresponding to $\pi_1(M_{i,j})$ are homeomorphic to $\text{int}(M_{i,j})$.
\elem
\bpf

We first claim that $M$ is atoroidal. Let $T\subset M$ be an essential torus. Since $T$ is compact it intersects at most finitely many $\set{S_i}_{i\in\Z}$. Moreover, up to an isotopy we can assume that $T$ is transverse to all $S_i$ and that it minimises $\abs{\pi_0(T\cap\cup_{i\in\Z} S_i)}$. If $T$ does not intersect any $S_i$ we have that it is contained in a submanifold homeomorphic to $Y$, see Remark \ref{fact0}, which is atoroidal and so $T$ isn't essential.

 Since both the $S_i$'s and $T$ are incompressible by our minimality condition we have that the components of the intersection $T\cap S_i$ are essential pairwise disjoint simple closed curves in $T,S_i$. Thus, $T$ is decomposed by $\cup_{i\in\Z} T\cap S_i$ into finitely many parallel annuli. Consider $S_k$ such that $T\cap S_k\neq \emp$ and $\forall n\geq k: T\cap S_n=\emp$. Then $T$ cannot intersect $S_k$ in only one component, so it has to come back through $S_k$. Thus, we have an annulus $A\subset T$ that has both boundaries in $S_k$ and is contained in a submanifold of $M$ homeomorphic to $Y$. The annulus $A$ gives an isotopy between isotopic curves in $\partial Y$ and is therefore boundary parallel, see Remark \ref{fact0}. Hence, by an isotopy of $T$ we can reduce $\abs{\pi_0(T\cap\cup_{i\in\Z} S_i)}$ contradicting the fact that it was minimal and non-zero. Therefore, $M$ is atoroidal.
 
 \vspace{0.3cm}
 
 \paragraph{Claim:} The $M_{i,j}$ are hyperbolizable.
 
 \vspace{0.3cm}
 \bpfc
Since $M$ is atoroidal and for $i<j$ the $M_{i,j}$ are $\pi_1$-injective submanifolds they are also atoroidal. Moreover, since the $M_{i,j}$ are compact manifolds with infinite $\pi_1$ they are hyperbolizable by Thurston's Hyperbolization Theorem \cite{Kap2001}.\epfc

\vspace{0.3cm}

The manifold $M$ is exhausted by the hyperbolizable $\pi_1$-injective submanifolds $M_i\eqdef M_{-i,i}$.
\vspace{0.3cm}

 \paragraph{Claim:} The manifold $M$ is locally hyperbolic.
 
 \vspace{0.3cm}
 \bpfc
To do so it suffices to show that given any finitely generated $H\subgroup \pi_1(M)$ the cover $M(H)$ corresponding to $H$ factors through a cover $N\twoheadrightarrow M$ that is hyperbolizable. Let $\gamma_1,\dotsc,\gamma_n\subset M$ be loops generating $H$. Since the $M_i$ exhaust $M$ we can find some $i\in\N$ such that $\set{\gamma_k}_{1\leq k\leq n} 	\subset M_i$, hence the cover corresponding to $H$ factors through the cover induced by $\pi_1(M_i)$. We now want to show that the cover $M (i)$ of $M$ corresponding to $\pi_1(M_i)$ is hyperbolizable.

Since $\pi:M\twoheadrightarrow X$ is an infinite cyclic cover of $X$ we have that $M(i)$ is the same as the cover of $X$ corresponding to $\pi_*(\pi_1(M_i))$. The resolution of the Tameness \cite{AG2004,CG2006} and the Geometrization conjecture \cite{Per2003.1,Per2003.2,Per2003.3} imply the Simon's conjecture\footnote{Final steps completed by Long and Reid, see \cite{Ca2010}.}, that is: covers of compact irreducible 3-manifolds with finitely generated fundamental groups are tame \cite{Ca2010,Si1976}. Therefore, since $X$ is compact by the Simon's Conjecture we have that $M(i)$ is tame. The submanifold $M_i\hookrightarrow M$ lifts homeomorphically to $\widetilde M_i\hookrightarrow M(i)$. By Whitehead's Theorem \cite{Ha2002} the inclusion is a homotopy equivalence, hence $\widetilde M_i$ forms a Scott core for $M(i)$. Thus, since $\partial \widetilde M_i$ is incompressible and $M(i)$ is tame we have that $M(i)\cong \text{int}(M_i)$ and so it is hyperbolizable. \epfc

Which concludes the proof. \epf

\brem\label{facts1}
Note that in the manifold $M$ the surfaces $S_i,S_j$ have no homotopic simple closed curve except for the loops $\gamma_i\eqdef S_i\cap A$. If not we would have an embedded cylinder $C$ not homotopic into $A_{i,j}$ which contradicts the fact that the characteristic submanifold of $M_{i,j}$ is given by a thickening of $A_{i,j}$.  In particular this gives us the important fact that for any homotopy equivalence $f:M\rar N$ and any essential subsurface $F\subset S_i$ not isotopic to a neighbourhood of $\gamma$ we cannot homotope $f(F)$ through any $f(S_j)$ for $i\neq j$.

\erem

\blem[Homotopy Equivalences]\label{homeq} Given a tame 3-manifold $N$ let $\overline N$ be its compactification and let $g: M_{i,j}\rar \overline N$ be a homotopy equivalence. Then, there exists a homotopy equivalence $f: M_{i,j}\rar\overline N$ such that $f\simeq g$ and:

\begin{enumerate}
\item $f(S_k)$ is embedded for all $i\leq k\leq j$;

\item there are essential subsurfaces $ T_m, T_n$ of $S_m, S_n$, respectively, whose components are homeomorphic to punctured tori, and where for all $i\leq m < k < n\leq j$, the images $f(T_n), f(T_m)$ are separated in $N$ by $f(S_k)$. Moreover, the same holds for any surface $\Sigma_k\iso f(S_k)$ intersecting $f(S_n),f(T_m)$ minimally.
\end{enumerate}
\elem
\bpf
By Lemma \ref{dehnflip} we get that $g:M_{i,j}\rar \overline N$ is either homotopic to a homeomorphism $f$ or is given by a Dehn flip along the annulus $A_{i,j}$ of $M_{i,j}$. If $N$ is homeomorphic to $M_{i,j}$ we have nothing to do since the required map $f$ is the homeomorphism and (1) and (2) are true for $M_{i,j}$.

Therefore, we only need to deal with the case in which $f: M_{i,j}\rar\overline N$ is a Dehn flip of $M$ along the annulus $A_{i,j}$. We will now explicitly write the Dehn flip  $f$. Let $V\cong S^1\times I_s\times I_t$ be a regular neighbourhood of the annulus $A_{i,j}$ in $M_{i,j}$ such that $V\cap S_k$, $i\leq k\leq j$, are regular neighbourhoods $S^1\times\set {s_k}\times I_t$ of $\gamma$ in $S_k$. Similarly let $W\cong S^1\times I\times I$ be a regular neighbourhood of $A_{i,j}$ in $\overline N$. Let $F:V\rar W$ be given by:

\be F(x,s,t)\eqdef \begin{cases} (x,2t(1-s)+(1-2t)s,t),\quad 0\leq t\leq \f 1 2\\ (x,(2-2t)(1-s)+(2t-1)s,t),\quad \f 12\leq t\leq 1\end{cases}\ee

and $f:M_{i,j}\rar\overline N$ be the homotopy equivalence obtained by extending $F$ via the homeomorphism of $\overline{M_{i,j}\setminus V}\rar \overline{\overline N\setminus W}$ coming from Lemma \ref{dehnflip}. Moreover, for $M_{i,j}'\eqdef \overline {M_{i,j}\setminus V}$ the homeomorphism $F$ is the identity on $V\cap M_{i,j}'$. Then $f$ realises the Dehn flip from $M_{i,j}$ to $\overline N$. The homeomorphism of $M'_{i,j}$ preserves the order of the surfaces, it is the identity on $\partial V\cap  M_{i,j}'$, hence for all $i\leq k\leq j$ the surfaces $f(S_k)$ are embedded. This concludes the proof of (1).

\begin{center}\begin{figure}[h!]
						\def\svgwidth{250pt}
						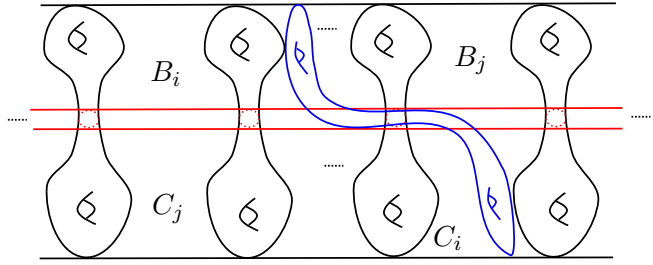
						\caption{The surface in blue is an embedded push-off of $f(S_i)$ in $N\simeq\text{int}( M_{i,j})$ when $\overline N$ is a Dehn flip of $M_{i,j}$.}	\label{fig2}
						\end{figure}
\end{center}

For (2) note that for all $i\leq k\leq j$ we have that $S_k\setminus V$ is given by two essential punctured tori $T_k^\pm$. Moreover, for all $i\leq n\neq k \leq j$ we have that the essential tori $T_n^\pm$ are separated by $f(S_k)$, this again follows from the fact that $f$ is the identity on $\partial V\cap  M_{i,j}'$ and so it preserves ordering. Thus, we always see, up to isotopy, the following configuration:
 
\begin{center}\begin{figure}[h!]
						\def\svgwidth{250pt}
						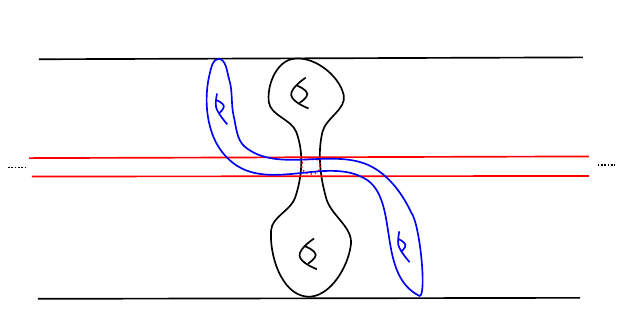
												\end{figure}
\end{center}

Finally if $\Sigma_k\iso f(S_k)$ and intersects $f(S_n), f(S_m)$ minimally we have that all components of intersections of $f(S_n), f(S_m)$ and $\Sigma_k$ are isotopic to the intersection of $f(S_n),f(S_m)$ with the annulus $A_{i,j}$. If the subsurfaces $f(T_n),f(T_m)$ are not separated by $\Sigma_k$ it means that they all lie in the same component $N_1$ of $N\setminus \Sigma_k$. 

By the isotopy extension Theorem \cite[8.1.3]{Hi1976} we have $\Sigma_n,\Sigma_m$ isotopic to $f(S_n)$ and $f(S_m)$ respectively and subsurfaces $T_n'$, $T_m'$ isotopic to $f(T_n),f(T_m)$ that are separated in $N$ by $\Sigma_k$. Therefore, we can find a simple closed loop $\alpha$ in one of the essential subsurfaces $T_n',T_m'$ that is not contained in $N_1$. Since we assumed that $f(T_n),f(T_m)$ are contained in $N_1$ the loop $\alpha$ is homotopic into $\Sigma_k$ contradicting Remark \ref{facts1}.\epf

%
%
%

		\bdefi	 Given a hyperbolic 3-manifold $M$, a  \textit{useful simplicial hyperbolic surface} is a surface $S$ with a 1-vertex triangulation $\mathcal T$, a preferred edge $e$ and a map $f:S\rar M$, such that:
			\begin{enumerate}
			 
			\item $f(e)$ is a geodesic in $M$; 
			\item every edge of $\mathcal T$ is mapped to a geodesic segment in $M$; 
			\item the restriction of $f$ to every face of $\mathcal T$ is a totally geodesic immersion.
			\end{enumerate}
				\edefi
By \cite{Ca1996,Bo1986} every $\pi_1$-injective map $f:S\rar M$ with a 1-vertex triangulation with a preferred edge can be homotoped so that it becomes a useful simplicial surface. Moreover, with the path metric induced by $M$ a useful simplicial surface is negatively curved and the map becomes $1$-Lipschitz.

\blem\label{lem2}
Let $N\cong\hyp\Gamma$ be a hyperbolic 3-manifold homotopy equivalent to $M$ then $\gamma$ is represented by a parabolic element in $\Gamma$.
\elem
\bpf Assume that $f_*(\gamma)$ is represented by a hyperbolic element and let $A\subset M$ be the essential bi-infinite annulus obtained as the limit of the $A_{-i,i}$. Since all $S_i$ are incompressible in $M$ and $f$ is a homotopy equivalence the maps: $f:S_i\rar N$ are $\pi_1$-injective. Let $\tau_i$ be 1-vertex triangulations of $S_i$ realising $\gamma_i\eqdef S_i\cap A$ as an edge in the 1-skeleton. Then, by \cite{Ca1996,Bo1986} we can realise the maps $f:S_i\rar N$ by useful simplicial hyperbolic surfaces $\Sigma_i\subset N$ such that $\Sigma_i\simeq f(S_i)$ and the image of $\gamma_i$ is the unique geodesic representative $\overline\gamma$ of $f_*(\gamma)$.

In the simplicial hyperbolic surfaces $\Sigma_i$ a maximal one-sided collar neighbourhood of $\overline\gamma$ has area bounded by the total area of $\Sigma_i$. Since the simplicial hyperbolic surfaces are all genus two by Gauss-Bonnet we have that $ A(\Sigma_i)\leq 2 \pi\abs{\chi(\Sigma_i)}=4\pi$. Therefore, the radius of a one-sided collar neighbourhood is uniformly bounded by some constant $K=K(\chi(\Sigma_i),\ell_N(\overline\gamma))<\infty$ for $\ell_N(\overline\gamma)$ the hyperbolic length of $\gamma$ in $N$. Then for $\xi>0$ in the simplicial hyperbolic surface $\Sigma_i$ the $K+\xi$ two sided neighbourhood of $\overline\gamma$ is not embedded and contains an essential 4-punctured sphere $X_i\subset\Sigma_i$. Since simplicial hyperbolic surfaces are $1$-Lipschitz the $4$-punctured sphere is contained in a $K+\xi$ neighbourhood $C$ of $\overline \gamma$ in $N$. The curves $\alpha_i$ obtained by joining the seams of the pants decomposition of $X_i$ induced by $\gamma_i$ have length bounded by $L\eqdef 2K+2\xi+\ell_N(\overline\gamma)$. Since there are infinitely many $\alpha_i$ and they are all homotopically distinct we have that $\Gamma$ is not discrete since the $\alpha_i$ move a lift of $\overline\gamma$ a uniformly bounded amount.

\epf

\bprop\label{prop3}
Let $f:M\rar N$ be a homotopy equivalence, then for all $i$ the maps: $f:S_i\rar N$ have embedded representatives $\Sigma_i\simeq f(S_i)$ in $N$. 
\eprop
\bpf Fix a triangulation $\tau$ of $N$. Since the $f(S_i)$ are $\pi_1$-injective by taking a refinement $\tau_i$ of $\tau$ outside a pre-compact neighbourhood $U_i$ of $f(S_i)$ we can homotope $f(S_i)$ to be a PL-least area surface $\Sigma_i$ with respect to a weight system induce by $\tau_i$, see \cite{FHS1983,Kap2001,JR1998}. We now want to show that they are embedded. To do so it suffices to show that they have embedded representatives in some cover, see \cite{FHS1983,Kap2001,JR1998}. Consider the cover $\pi_i:N_{i,i+1}\twoheadrightarrow N$ with the triangulation $\tilde\tau_i$ induced by $\tau_i$. The PL-least area surface $\Sigma_i$ lifts homeomorphically to $\tilde\Sigma_i\subset N_{i,i+1}$ and it is still minimal with respect to $\tilde\tau_i$. Since $f(S_i)\simeq\Sigma_i$ has embedded representatives in $N_{i,i+1}$, see Lemma \ref{homeq}, by \cite{FHS1983,Kap2001,JR1998} we have that $\tilde\Sigma_i$ is embedded as well, hence $\Sigma_i=\pi_i(\tilde\Sigma_i)$ is.
\epf

We will now prove Theorem \ref{maintheorem} which we now restate.

\bthm
The manifold $M $ is locally hyperbolic and without divisible elements in $\pi_1(M)$ but is not homotopy equivalent to any hyperbolic 3-manifold.
\ethm
\bpf By Lemma \ref{nodivis} and Lemma \ref{lochyp} we only need to show that $M$ is not homotopy equivalent to any hyperbolic 3-manifold $N$. The proof will be by contradiction. Assume that we have a homotopy equivalence $f:M\rar N$ for $N\cong\hyp\Gamma$ a hyperbolic 3-manifold. By Lemma \ref{lem2} we have that for $\gamma$ the element of $\pi_1(M)$ generating the fundamental group of the bi-infinite essential annulus $A\subset M$ $f_*(\gamma)$ is represented by a parabolic element in $\Gamma$ (with an abuse of notation we will refer to this element by $\gamma$ as well). Thus, in $N$ we have a cusp $E\eqdef E_\gamma$ corresponding to $\gamma$. Moreover, Proposition \ref{prop3} gives us a collection $\set{\Sigma_i}_{i\in\Z}$ of embedded genus two surfaces contained in neighbourhoods $U_i$ of $f(S_i)$ in $N$ such that $\Sigma_i\simeq f(S_i)$. Moreover the $\Sigma_i$'s are incompressible and separating in $N$. The fact that they are separating follows from $f_*$ being an isomorphism in homology and the fact that the $S_i$'s are not dual to 1-cycles in $M$, similarly they are incompressible since the $S_i$ are and $f$ is a homotopy equivalence. Therefore, if we take the surface $\Sigma_0$ we have that $N\vert \Sigma_0\eqdef N\setminus \text{int}(N_r(\Sigma_0))$ is given by two manifolds $N_1,N_2$ with boundary a surface isotopic in $N$ to $\Sigma_0$.

The element $\gamma\in\Gamma$ is parabolic with cusp $E$ and each $\Sigma_i$ has a simple closed loop $\gamma_i$ homotopic in $N$ to $\gamma$ such that $\overline{\Sigma_i\setminus N_r(\gamma_i)}$ is given by two punctured tori $T^\pm_i$ with boundary isotopic to $\gamma_i$.  Without loss of generality we can assume that $E\subset N_2$. Moreover, up to an isotopy of each $\Sigma_i$ we can also assume that for all $i\in\Z$ the surfaces $\Sigma_i$ are transverse to $\Sigma_0$ and that $\abs{\pi_0(\Sigma_i\cap\Sigma_0)}$ is minimal.

\vspace{0.3cm}

\paragraph{Claim:} Every component $\alpha\in\pi_0(\Sigma_i\cap\Sigma_0)$ is isotopic to $\gamma_i$ and $\gamma_0$.

\vspace{0.3cm}

\bpfc  Since $\Sigma_i,\Sigma_0$ are incompressible and we minimised $\Sigma_i\cap\Sigma_0$ we have that every $\alpha\in\pi_0(\Sigma_i\cap\Sigma_0)$ has to be essential in both surfaces. By Remark \ref{facts1} we have that the only simple closed curve in $S_i$ homotopic into $S_0$ is $\gamma_i$ which is homotopic to $\gamma_0$. \epfc

Thus in $N$ we have that the punctured tori $T_i^\pm$ to $\Sigma_i$ are either on the same side of $\Sigma_0$ or on opposite sides:

\begin{center}\begin{figure}[h!]
						\def\svgwidth{300pt}
						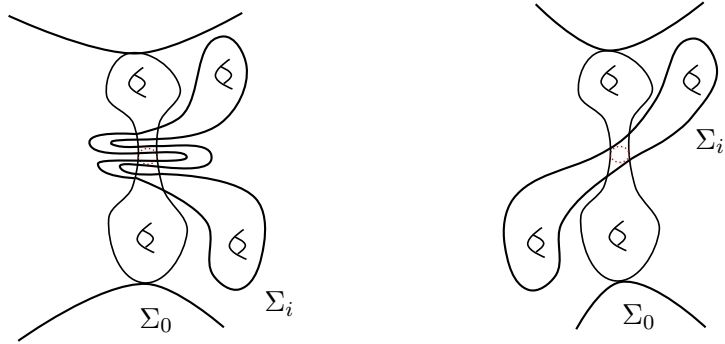
						\caption{Possible configurations of the embedded surfaces $\Sigma_0,\Sigma_i$ in $N$.}	\label{fig2}
						\end{figure}
\end{center}

Moreover, all the components of $\Sigma_i\cap\Sigma_0$ are contained in neighbourhoods of $\gamma_0$ and $\gamma_i$.

\vspace{0.3cm}

\paragraph{Claim:} There are infinitely many punctured tori $\set{T_n}_{n\in\N}\subset N_1$ such that $T_n$ is a component of $T_{i_n}^\pm$.

\vspace{0.3cm}

\bpfc Consider $\Sigma_{-i},\Sigma_0,\Sigma_{i}$ and a cover $\pi_j:N_{-i,i}\twoheadrightarrow N$ corresponding to the subgroup $f_*(\pi_1(M_{-j,j}))\subset\pi_1(N)$ where $\Sigma_0\cup\Sigma_i\cup\Sigma_{-i}$ lifts homeomorphically. Assume that there only finitely many $T_i^\pm$ that are contained in $N_1$. Then, for infinitely many $T_i^\pm$ in the covers $N_j$ we see the following configuration:

\begin{center}\begin{figure}[h!]
						\def\svgwidth{300pt}
						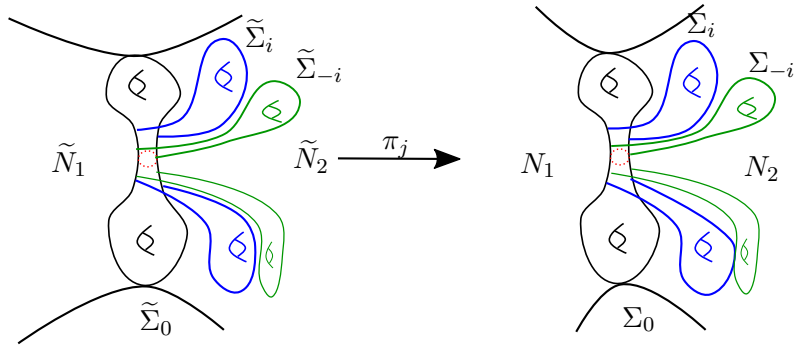
						\caption{The tori $T_i^\pm$ are marked by the surfaces $\Sigma_i$ that they are subsurfaces of.}	\label{fig2}
						\end{figure}
\end{center}

Let $g$ be as in Lemma \ref{homeq} and homotopic to the homotopy equivalence $\tilde f: M_{-j,j}\rar \overline N_{-j,j}$. Since $\tilde \Sigma_k$ and $g(S_k)$ are incompressible closed surfaces by \cite{Wa1968} we have that $\tilde\Sigma_k\simeq g(S_k)$ and by (2) of Lemma \ref{homeq} we have that $\tilde\Sigma_k$ separates the punctured tori $g(Q_i^\pm),g(Q_{-i}^\pm)$ in $g(S_i)$ and $g(S_{-i})$. Thus, we have a punctured torus, say $f(Q_i^+)$, that is contained in $\tilde N_1$ and such that the corresponding punctured torus $\tilde T_i^+\subset\tilde\Sigma_i$ is contained in $\tilde N_2$. Let $\alpha\subset T_i^+$ be any essential non peripheral curve, then since $\alpha$ is homotopic into $f(Q_i^+)$ and $\tilde\Sigma_0$ separates we have that $\alpha$ is homotopic into $\tilde\Sigma_0$ contradicting Remark \ref{facts1}. Therefore, we have infinitely many punctured tori $\set{T_n}_{n\in\N}$ with boundary $\partial T_n\eqdef \gamma_n$ isotopic to $\gamma$ such that $\forall n: T_n\subset N_1$. \epfc

We can now reach a contradiction with the fact that $\Gamma$ is a discrete group. Let $\mu\eqdef\text{inj}_N(\Sigma_0)$ and let $\epsilon\eqdef \min\set{\mu,\mu_3}$ where $\mu_3$ is the 3-dimensional Margulis constant (see \cite{BP1992}). Since the $T_n$'s are $\pi_1$-injective by picking a $1$-vertex triangulation $\tau_n$ with preferred edge corresponding to $\gamma_n$ we can realise the $T_n$'s by useful simplicial hyperbolic surfaces $P_n$. 

The surfaces $P_n$ are mapping $\gamma_n$ to the cusps $E$ and cannot be homotoped through $\Sigma_0$ since again we would contradict Remark \ref{facts1}. Hence, we get that for all $n$ $P_n\cap \Sigma_0\neq\emp$. Let $x_n\in P_n\cap\Sigma_n$, then $\text{inj}_{x_n}(P_n)\geq\epsilon$. Since the surfaces $P_n$ are negatively curved by the Bounded Diameter Lemma \cite{Bo1986,Th1978}  we get that we can find loops $\alpha_n,\beta_n\in\pi_1(P_n,x_n)$ whose length is bounded by $\f 8 \epsilon$ and such that they generate a rank two free group. Since $\langle \alpha_n,\beta_n\rangle\cong \mathbb F_2$ we have that at least one of $\alpha_n,\beta_n$ is not homotopic to $\gamma_n$. Without loss of generality we can assume that this element is always $\alpha_n$. By Remark \ref{facts1} the collection $\set{\alpha_n}_{n\in\N}$ are all distinct elements in $\Gamma$. Moreover, we have that $\ell_N(\alpha_n)\leq\f 8\epsilon$. Let $D\eqdef\text{diam}(\Sigma_0)$, pick $x\in\Sigma_0$ and fix $\tilde x$ to be a lift in $\tilde \Sigma_0\subset\bH$. Then for lifts $\tilde x_n$ of $x_n$ we have that:

\begin{align*}
d_{\bH}(\alpha_n(\tilde x),\tilde x)&\leq d_{\bH}(\alpha_n(\tilde x),\alpha_n(\tilde x_n))+d_{\bH}(\alpha_n(\tilde x_n),\tilde x_n)+d_{\bH}(\tilde x_n,\tilde x)\\
&\leq D+\f 8 \epsilon +D\\
&=2D+\f 8 \epsilon
\end{align*} 

Thus the family $\set{\alpha_n}_{n\in\N}$ has an accumulation point in $\text{PSL}_2(\mathbb C)$ contradicting the discreetness of $\Gamma$. \epf

\thispagestyle{empty}
{\small
\markboth{References}{References}
\bibliographystyle{alpha}
\bibliography{mybib}{}
}

\noindent Department of Mathematics, Boston College.

\noindent 140 Commonwealth Avenue Chestnut Hill, MA 02467.

\noindent Maloney Hall
\newline \noindent
email: \texttt{cremasch@bc.edu}
	
	\end{document}

%% file: packages.tex
\usepackage{latexsym,enumitem}
\usepackage{amssymb}
\usepackage[cp850]{inputenc}
\usepackage{epsfig}
\usepackage{psfrag}
\usepackage{amsthm}
\usepackage{amscd}
\usepackage{amsmath}
\usepackage{amsfonts}
\usepackage{graphics,caption}
\usepackage[all]{xy}
\usepackage{etoolbox}
\patchcmd{\quote}{\rightmargin}{\leftmargin 2em \rightmargin}{}{}
 \usepackage[sort,numbers]{natbib}

\usepackage[english]{babel}
\usepackage{comment}
\usepackage{color}
\usepackage[colorlinks]{hyperref}

\let\phi\varphi

\DeclareMathOperator{\eqdef}{\doteq\,} 
\newcommand{\f}[2]{\frac{#1}{#2}} 

\let\epsilon\varepsilon
\let\subset\subseteq
\newcommand{\be}{\begin{equation*}}
 \newcommand{\ee}{\end{equation*}}
\newcommand{\bpf}{\begin{dimo}}
\newcommand{\epf}{\end{dimo}}
\newcommand{\bdefi}{\begin{defin}}
\newcommand{\edefi}{\end{defin}}
\newcommand{\bthm}{\begin{thm}}
\newcommand{\ethm}{\end{thm}}
\newcommand{\blem}{\begin{lem}}
\newcommand{\elem}{\end{lem}}
\newcommand{\bcor}{\begin{cor}}
\newcommand{\ecor}{\end{cor}}

\newcommand{\bprop}{\begin{prop}}
\newcommand{\eprop}{\end{prop}}
\newcommand{\bese}{\begin{ese}}
\newcommand{\eese}{\end{ese}}
\newcommand{\brem}{\begin{rem}}
\newcommand{\erem}{\end{rem}}
\newcommand{\bpfc}{\begin{dimoclaim}}
\newcommand{\epfc}{\end{dimoclaim}}

\newcommand{\rar}{\rightarrow} 
\newcommand{\imm}{\looparrowright}

\newcommand{\abs}[1]{\left\lvert#1\right\rvert}						
\newcommand{\set}[1]{\left\{#1\right\}}					
	
\newcommand{\quotient}[2]{\left.\raisebox{.1em}{$#1\!$}\middle/\raisebox{-.1em}{$#2$}\right.}

\DeclareMathOperator{\emp}{\varnothing} 
\DeclareMathOperator{\N}{\mathbb N}			
\DeclareMathOperator{\Z}{\mathbb Z}			

\newcommand{\diffeo}{\xrightarrow{{}_{\simeq}}} 
\newcommand{\bH}{\mathbb H^3}
 \newcommand{\hyp}[1]{\quotient{\bH}{#1}}
  \newcommand{\iso}{\overset{\text{iso}}\simeq}
\newcommand{\subgroup}{\leqslant}

\newenvironment{quot}
{
	\vspace{-0.2cm}
	\vspace{0.2cm}
}

\theoremstyle{definition}
\newtheorem{d1}{Definition}[section] 

\newenvironment{defin}
{
	\begin{quot}
		\begin{d1}
		}
		{\end{d1}
	\end{quot}

}

\theoremstyle{definition}
\newtheorem{r1}[d1]{Remark}

\newenvironment{rem}
{
	\begin{quot}
		\begin{r1}
		}
		{\end{r1}
	\end{quot}
}

\theoremstyle{definition}
\newtheorem{e1}[d1]{Exercise}

\theoremstyle{definition}
\newtheorem{ese1}[d1]{Example}

\newenvironment{ese}
{
	\begin{quot}
		\begin{ese1}
	}
	{	
		\end{ese1}
	\end{quot}
}

\theoremstyle{definition}

\theoremstyle{definition}
\newtheorem{f2}[d1]{Fact}

\theoremstyle{definition}

\theoremstyle{definition}

\theoremstyle{definition}
\newtheorem{t1}[d1]{Theorem}

\newenvironment{thm}
{
	\begin{quot}
		\begin{t1}}
		{\end{t1}
	\end{quot}
}

\theoremstyle{definition}
\newtheorem*{T1*}{Theorem}

\newenvironment{teor*}
{
	\begin{quot}
		\begin{T1*}}
		{\end{T1*}
	\end{quot}
}

\newenvironment{dimo}
{\begin{proof}[Proof]
	}
	{\end{proof}}

\newenvironment{dimoclaim}{\emph{Proof of Claim:}\;}{\hfill$\square$}

	\theoremstyle{definition}
	\newtheorem{l1}[d1]{Lemma}
	
	\newenvironment{lem}
	{
		\begin{quot}
			\begin{l1}}
			{\end{l1}
		\end{quot}
	}
	\theoremstyle{definition}
	\newtheorem{p1}[d1]{Proposition}
	
	\newenvironment{prop}
	{
		\begin{quot}
			\begin{p1}}
			{\end{p1}
		\end{quot}
	}
	
	\theoremstyle{definition}
	\newtheorem{c1}[d1]{Corollary}
	
	\newenvironment{cor}
	{
		\begin{quot}
			\begin{c1}}
			{\end{c1}
		\end{quot}
	}

 \newtheorem*{Theorem*}{Theorem}
 \newtheorem*{Proposition*}{Proposition}
 \newtheorem*{Lemma*}{Lemma}

%% file: cover.pdf_tex
\begingroup%
  \makeatletter%
  \providecommand\color[2][]{%
    \errmessage{(Inkscape) Color is used for the text in Inkscape, but the package 'color.sty' is not loaded}%
    \renewcommand\color[2][]{}%
  }%
  \providecommand\transparent[1]{%
    \errmessage{(Inkscape) Transparency is used (non-zero) for the text in Inkscape, but the package 'transparent.sty' is not loaded}%
    \renewcommand\transparent[1]{}%
  }%
  \providecommand\rotatebox[2]{#2}%
  \ifx\svgwidth\undefined%
    \setlength{\unitlength}{189.37174988bp}%
    \ifx\svgscale\undefined%
      \relax%
    \else%
      \setlength{\unitlength}{\unitlength * \real{\svgscale}}%
    \fi%
  \else%
    \setlength{\unitlength}{\svgwidth}%
  \fi%
  \global\let\svgwidth\undefined%
  \global\let\svgscale\undefined%
  \makeatother%
  \begin{picture}(1,0.62009657)%
    \put(0,0){\includegraphics[width=\unitlength,page=1]{cover.pdf}}%
    \put(0.70113974,0.26026686){\color[rgb]{0,0,0}\makebox(0,0)[lb]{\smash{$S$}}}%
    \put(0.45175018,0.40471116){\color[rgb]{0,0,0}\makebox(0,0)[lb]{\smash{$B$}}}%
    \put(0.81987104,0.45319396){\color[rgb]{0,0,0}\makebox(0,0)[lb]{\smash{$C$}}}%
    \put(0.67455548,0.40953018){\color[rgb]{0,0,0}\makebox(0,0)[lb]{\smash{$T$}}}%
    \put(0,0){\includegraphics[width=\unitlength,page=2]{cover.pdf}}%
  \end{picture}%
\endgroup%

%% file: fig1.pdf_tex
\begingroup%
  \makeatletter%
  \providecommand\color[2][]{%
    \errmessage{(Inkscape) Color is used for the text in Inkscape, but the package 'color.sty' is not loaded}%
    \renewcommand\color[2][]{}%
  }%
  \providecommand\transparent[1]{%
    \errmessage{(Inkscape) Transparency is used (non-zero) for the text in Inkscape, but the package 'transparent.sty' is not loaded}%
    \renewcommand\transparent[1]{}%
  }%
  \providecommand\rotatebox[2]{#2}%
  \ifx\svgwidth\undefined%
    \setlength{\unitlength}{300bp}%
    \ifx\svgscale\undefined%
      \relax%
    \else%
      \setlength{\unitlength}{\unitlength * \real{\svgscale}}%
    \fi%
  \else%
    \setlength{\unitlength}{\svgwidth}%
  \fi%
  \global\let\svgwidth\undefined%
  \global\let\svgscale\undefined%
  \makeatother%
  \begin{picture}(1,0.5)%
    \put(0,0){\includegraphics[width=\unitlength,page=1]{fig1.pdf}}%
    \put(0.12818451,0.00770095){\color[rgb]{0,0,0}\makebox(0,0)[lb]{\smash{$S_{i-1}$}}}%
    \put(0.37369137,0.00412758){\color[rgb]{0,0,0}\makebox(0,0)[lb]{\smash{$S_i$}}}%
    \put(0.59469987,0.00318023){\color[rgb]{0,0,0}\makebox(0,0)[lb]{\smash{$S_{i+1}$}}}%
    \put(0.83400718,0.00575224){\color[rgb]{0,0,0}\makebox(0,0)[lb]{\smash{$S_{i+2}$}}}%
    \put(0.27693575,-0.0518138){\color[rgb]{0,0,0}\makebox(0,0)[lb]{\smash{}}}%
    \put(0.96563314,0.24770489){\color[rgb]{1,0,0}\makebox(0,0)[lb]{\smash{$A$}}}%
    \put(0,0){\includegraphics[width=\unitlength,page=2]{fig1.pdf}}%
    \put(0.44709889,0.34280405){\color[rgb]{0,0,0}\makebox(0,0)[lb]{\smash{$B_i$}}}%
    \put(0.46046823,0.12121303){\color[rgb]{0,0,0}\makebox(0,0)[lb]{\smash{$C_i$}}}%
    \put(0.68976406,0.33813441){\color[rgb]{0,0,0}\makebox(0,0)[lb]{\smash{$B_{i+1}$}}}%
    \put(0.69470245,0.14129099){\color[rgb]{0,0,0}\makebox(0,0)[lb]{\smash{$C_{i+1}$}}}%
    \put(0.20170702,0.3392324){\color[rgb]{0,0,0}\makebox(0,0)[lb]{\smash{$B_{i-1}$}}}%
    \put(0.22266653,0.1330647){\color[rgb]{0,0,0}\makebox(0,0)[lb]{\smash{$C_{i-1}$}}}%
    \put(0.4605279,0.45685764){\color[rgb]{0,0,0}\makebox(0,0)[lb]{\smash{$Y_i$}}}%
    \put(0.25633382,0.46319904){\color[rgb]{0,0,0}\makebox(0,0)[lb]{\smash{$Y_{i-1}$}}}%
    \put(0.69544321,0.45929198){\color[rgb]{0,0,0}\makebox(0,0)[lb]{\smash{$Y_{i+1}$}}}%
  \end{picture}%
\endgroup%

%% file: fig2.pdf_tex
\begingroup%
  \makeatletter%
  \providecommand\color[2][]{%
    \errmessage{(Inkscape) Color is used for the text in Inkscape, but the package 'color.sty' is not loaded}%
    \renewcommand\color[2][]{}%
  }%
  \providecommand\transparent[1]{%
    \errmessage{(Inkscape) Transparency is used (non-zero) for the text in Inkscape, but the package 'transparent.sty' is not loaded}%
    \renewcommand\transparent[1]{}%
  }%
  \providecommand\rotatebox[2]{#2}%
  \ifx\svgwidth\undefined%
    \setlength{\unitlength}{300bp}%
    \ifx\svgscale\undefined%
      \relax%
    \else%
      \setlength{\unitlength}{\unitlength * \real{\svgscale}}%
    \fi%
  \else%
    \setlength{\unitlength}{\svgwidth}%
  \fi%
  \global\let\svgwidth\undefined%
  \global\let\svgscale\undefined%
  \makeatother%
  \begin{picture}(1,0.5)%
    \put(0,0){\includegraphics[width=\unitlength,page=1]{fig2.pdf}}%
   \put(0,0){\includegraphics[width=\unitlength,page=2]{fig2.pdf}}%
    \put(0.22838224,0.29026245){\color[rgb]{0,0,0}\makebox(0,0)[lb]{\smash{$B_i$}}}%
    \put(0.687955,0.31344778){\color[rgb]{0,0,0}\makebox(0,0)[lb]{\smash{$B_{j}$}}}%
    \put(0,0){\includegraphics[width=\unitlength,page=3]{fig2.pdf}}%
    \put(0.23067567,0.08877578){\color[rgb]{0,0,0}\makebox(0,0)[lb]{\smash{$C_{j}$}}}%
    \put(0.65609227,0.04132797){\color[rgb]{0,0,0}\makebox(0,0)[lb]{\smash{$C_i$}}}%
    \put(0,0){\includegraphics[width=\unitlength,page=4]{fig2.pdf}}%
  \end{picture}%
\endgroup%

%% file: fig2_5.pdf_tex
\begingroup%
  \makeatletter%
  \providecommand\color[2][]{%
    \errmessage{(Inkscape) Color is used for the text in Inkscape, but the package 'color.sty' is not loaded}%
    \renewcommand\color[2][]{}%
  }%
  \providecommand\transparent[1]{%
    \errmessage{(Inkscape) Transparency is used (non-zero) for the text in Inkscape, but the package 'transparent.sty' is not loaded}%
    \renewcommand\transparent[1]{}%
  }%
  \providecommand\rotatebox[2]{#2}%
  \ifx\svgwidth\undefined%
    \setlength{\unitlength}{300bp}%
    \ifx\svgscale\undefined%
      \relax%
    \else%
      \setlength{\unitlength}{\unitlength * \real{\svgscale}}%
    \fi%
  \else%
    \setlength{\unitlength}{\svgwidth}%
  \fi%
  \global\let\svgwidth\undefined%
  \global\let\svgscale\undefined%
  \makeatother%
  \begin{picture}(1,0.5)%
    \put(0,0){\includegraphics[width=\unitlength,page=1]{fig2_5.pdf}}%
    \put(0.27693575,-0.0518138){\color[rgb]{0,0,0}\makebox(0,0)[lb]{\smash{}}}%
    \put(0.70018046,0.05932112){\color[rgb]{0,0,0}\makebox(0,0)[lb]{\smash{$f(S_n)$}}}%
    \put(0.46459005,0.43191253){\color[rgb]{0,0,0}\makebox(0,0)[lb]{\smash{$f(S_k)$}}}%
  \end{picture}%
\endgroup%

%% file: fig3.pdf_tex
\begingroup%
  \makeatletter%
  \providecommand\color[2][]{%
    \errmessage{(Inkscape) Color is used for the text in Inkscape, but the package 'color.sty' is not loaded}%
    \renewcommand\color[2][]{}%
  }%
  \providecommand\transparent[1]{%
    \errmessage{(Inkscape) Transparency is used (non-zero) for the text in Inkscape, but the package 'transparent.sty' is not loaded}%
    \renewcommand\transparent[1]{}%
  }%
  \providecommand\rotatebox[2]{#2}%
  \ifx\svgwidth\undefined%
    \setlength{\unitlength}{300bp}%
    \ifx\svgscale\undefined%
      \relax%
    \else%
      \setlength{\unitlength}{\unitlength * \real{\svgscale}}%
    \fi%
  \else%
    \setlength{\unitlength}{\svgwidth}%
  \fi%
  \global\let\svgwidth\undefined%
  \global\let\svgscale\undefined%
  \makeatother%
  \begin{picture}(1,0.5)%
    \put(0,0){\includegraphics[width=\unitlength,page=1]{fig3.pdf}}%
    \put(0.27693575,-0.0518138){\color[rgb]{0,0,0}\makebox(0,0)[lb]{\smash{}}}%
    \put(0,0){\includegraphics[width=\unitlength,page=2]{fig3.pdf}}%
    \put(0.19558563,0.0482404){\color[rgb]{0,0,0}\makebox(0,0)[lb]{\smash{$\Sigma_0$}}}%
    \put(0.80408241,0.05254377){\color[rgb]{0,0,0}\makebox(0,0)[lb]{\smash{$\Sigma_0$}}}%
    \put(0.35308895,0.06889658){\color[rgb]{0,0,0}\makebox(0,0)[lb]{\smash{$\Sigma_i$}}}%
    \put(0.89617459,0.27373709){\color[rgb]{0,0,0}\makebox(0,0)[lb]{\smash{$\Sigma_i$}}}%
     \end{picture}%
\endgroup%

%% file: fig4.pdf_tex
\begingroup%
  \makeatletter%
  \providecommand\color[2][]{%
    \errmessage{(Inkscape) Color is used for the text in Inkscape, but the package 'color.sty' is not loaded}%
    \renewcommand\color[2][]{}%
  }%
  \providecommand\transparent[1]{%
    \errmessage{(Inkscape) Transparency is used (non-zero) for the text in Inkscape, but the package 'transparent.sty' is not loaded}%
    \renewcommand\transparent[1]{}%
  }%
  \providecommand\rotatebox[2]{#2}%
  \ifx\svgwidth\undefined%
    \setlength{\unitlength}{300bp}%
    \ifx\svgscale\undefined%
      \relax%
    \else%
      \setlength{\unitlength}{\unitlength * \real{\svgscale}}%
    \fi%
  \else%
    \setlength{\unitlength}{\svgwidth}%
  \fi%
  \global\let\svgwidth\undefined%
  \global\let\svgscale\undefined%
  \makeatother%
  \begin{picture}(1,0.5)%
    \put(0,0){\includegraphics[width=\unitlength,page=1]{fig4.pdf}}%
    \put(0.27693575,-0.0518138){\color[rgb]{0,0,0}\makebox(0,0)[lb]{\smash{}}}%
    \put(0,0){\includegraphics[width=\unitlength,page=2]{fig4.pdf}}%
    \put(0.19558563,0.0482404){\color[rgb]{0,0,0}\makebox(0,0)[lb]{\smash{$\tilde\Sigma_0$}}}%
    \put(0.80408241,0.05254377){\color[rgb]{0,0,0}\makebox(0,0)[lb]{\smash{$\Sigma_0$}}}%
    \put(0.3282082,0.40506281){\color[rgb]{0,0,0}\makebox(0,0)[lb]{\smash{$\tilde\Sigma_i$}}}%
    \put(0.88511648,0.42357434){\color[rgb]{0,0,0}\makebox(0,0)[lb]{\smash{$\Sigma_i$}}}%
    \put(0,0){\includegraphics[width=\unitlength,page=3]{fig4.pdf}}%
    \put(0.50146637,0.27634579){\color[rgb]{0,0,0}\makebox(0,0)[lb]{\smash{$\pi_j$}}}%
    \put(0,0){\includegraphics[width=\unitlength,page=4]{fig4.pdf}}%
    \put(0.38945524,0.25276912){\color[rgb]{0,0,0}\makebox(0,0)[lb]{\smash{$\tilde N_2$}}}%
    \put(0.95765724,0.2422944){\color[rgb]{0,0,0}\makebox(0,0)[lb]{\smash{$N_2$}}}%
    \put(0,0){\includegraphics[width=\unitlength,page=5]{fig4.pdf}}%
    \put(0.96307603,0.37381288){\color[rgb]{0,0,0}\makebox(0,0)[lb]{\smash{$\Sigma_{-i}$}}}%
    \put(0.39234519,0.3669124){\color[rgb]{0,0,0}\makebox(0,0)[lb]{\smash{$\tilde\Sigma_{-i}$}}}%
    \put(0.08425167,0.24834588){\color[rgb]{0,0,0}\makebox(0,0)[lb]{\smash{$\tilde N_1$}}}%
    \put(0.67567568,0.24671765){\color[rgb]{0,0,0}\makebox(0,0)[lb]{\smash{$N_1$}}}%
  \end{picture}%
\endgroup%